\documentclass[10pt]{amsart}

\usepackage[T1]{fontenc}
\usepackage[ngerman,english]{babel}
\usepackage{amsmath}
\usepackage{amsthm}
\usepackage{amssymb}
\usepackage{amsfonts}
\usepackage{amsxtra}
\usepackage{color}
\usepackage[breaklinks]{hyperref}
\usepackage{enumerate}

\numberwithin{equation}{section}

\newcommand{\td}{\,\mathrm{d}}

\renewcommand\Re{\operatorname{Re}}
\renewcommand\Im{\operatorname{Im}}
\renewcommand{\det}{\textup{det}}

\newcommand{\rk}{\textup{rk}}

\newcommand{\GL}{\textup{GL}}

\newcommand{\SL}{\textup{SL}}
\newcommand{\SLpm}{\textup{S}^*\textup{L}}

\newcommand{\Sp}{\textup{Sp}}

\newcommand{\Mp}{\textup{Mp}}
\newcommand{\upO}{\textup{O}}
\newcommand{\upU}{\textup{U}}
\newcommand{\SO}{\textup{SO}}

\newcommand{\SU}{\textup{SU}}

\newcommand{\RR}{\mathbb{R}}

\newcommand{\CC}{\mathbb{C}}
\newcommand{\ZZ}{\mathbb{Z}}
\newcommand{\NN}{\mathbb{N}}

\renewcommand{\SS}{\mathbb{S}}
\newcommand{\HH}{\mathbb{H}}

\newcommand{\PP}{\mathbb{P}\hspace{0.4mm}}

\newcommand{\1}{\mathbf{1}}
\newcommand{\Ind}{\textup{Ind}}
\newcommand{\sgn}{\textup{sgn}}

\newcommand{\calO}{\mathcal{O}}

\newcommand{\calH}{\mathcal{H}}

\newcommand{\calL}{\mathcal{L}}

\newcommand{\frakg}{\mathfrak{g}}

\newcommand{\fraka}{\mathfrak{a}}

\newcommand{\frakh}{\mathfrak{h}}

\newcommand{\diag}{\textup{diag}}

\renewcommand{\max}{\textup{max}}

\newcommand{\ds}{\textup{ds}}

\theoremstyle{plain}
\newtheorem{theorem}{Theorem}

\theoremstyle{definition}

\theoremstyle{remark}
\renewenvironment{pf}{\begin{proof}[\sc Proof]}{\end{proof}}

\numberwithin{equation}{section}

\begin{document}

\title[Branching for small representations of $\GL(n,\CC)$]{Branching laws for small unitary representations of $\GL(n,\CC)$}

\author{Jan M\"ollers}
\author{Benjamin Schwarz}

\address{Institut for Matematiske Fag, Aarhus Universitet, Ny Munkegade 118, 8000 Aarhus C, Denmark}
\email{moellers@imf.au.dk}

\address{Institut f\"ur Mathematik, Universit\"at Paderborn, Warburger Str. 100, 33098 Paderborn, Germany}
\email{bschwarz@math.upb.de}
\begin{abstract}
The unitary principal series representations of $G=\GL(n,\CC)$ induced from a character of the maximal parabolic subgroup $P=(\GL(1,\CC)\times\GL(n-1,\CC))\ltimes\CC^{n-1}$ attain the minimal Gelfand--Kirillov dimension among all infinite-dimensional unitary representations of $G$. We find the explicit branching laws for the restriction of these representations to symmetric subgroups of $G$.
\end{abstract}



\maketitle

\section*{Introduction}

Branching laws describe the decomposition of a given irreducible unitary representation of a group $G$ into irreducible representations for a subgroup $H$. For reductive groups, say, every such representation can be decomposed into a direct integral of irreducible unitary representations for $H$. The explicit decomposition, however, is in general hard to determine. Further, the branching laws might have bad features such as infinite multiplicities, even if $(G,H)$ is a symmetric pair of reductive groups (see e.g. \cite[Example 5.5]{Kob00}). Therefore one has to single out certain subclasses of groups and representations in order to obtain nice branching laws.

As advocated by Kobayashi~\cite{Kob10} ``small representations'' contain ``large symmetries'' in their geometric realizations. This principle should be reflected in the branching; the restriction of a small representation of a non-compact group $G$ to a large subgroup $H$ is expected to have a simple branching law. Here the size of a unitary representation of a reductive group $G$ can be measured in terms of its Gelfand--Kirillov dimension. If the complexified Lie algebra $\frakg_\CC$ of $G$ does not contain a simple factor of type $A$ there are only finitely many irreducible unitary representations attaining the minimal Gelfand--Kirillov dimension among all non-trivial unitary representations \cite{GS05}. For the metaplectic group $G=\Mp(n,\RR)$ this is achieved by the well-known Segal--Shale--Weil representation; branching laws for this representation are well-studied, in particular in connection with Howe dual pairs, see e.g. \cite{How89,KV78}. More examples of branching laws for small representations are \cite{Dvo07,DS99,KO03,Li00,MS10,MO13,Mor08,PZ04,Sep07b,Sep07a,Zha01}.

On the other hand, for the type $A$ group $G=\GL(n,\RR)$ there exists a one-parameter family of irreducible unitary representations attaining the minimal Gelfand--Kirillov dimension. For these representations the explicit branching laws with respect to all symmetric subgroups $H\subseteq G$ were given by Kobayashi--{\O}rsted--Pevzner~\cite{KOP11}. In particular they investigate the branching to the subgroup $H=\Sp(m,\RR)$, $n=2m$, in connection with geometric analysis on various realizations of the representation. The corresponding geometric analysis for the group $H=\Sp(m,\CC)$ was carried out by Clare~\cite{Cla12}. The representations Clare studies occur as restrictions of irreducible unitary representations of the type $A$ group $G=\GL(n,\CC)$ that attain the minimal Gelfand--Kirillov dimension. In this paper we complete the picture by determining the explicit branching laws for the restriction of these small representations of $G=\GL(n,\CC)$ to arbitrary symmetric subgroups $H\subseteq G$.\\

Consider the unitary principal series representations $\pi_{i\lambda,k}^G$ of $G=\GL(n,\CC)$ induced from a character $\chi_{i\lambda,k}$, $\lambda\in\RR$, $k\in\ZZ$, of the parabolic subgroup $P=(\GL(1,\CC)\times\GL(n-1,\CC))\ltimes\CC^{n-1}$ (see Section~\ref{sec:UnitaryPrincipalSeriesRepresentations} for details). The arguments of \cite{MS12} show that the representations $\pi_{i\lambda,m}^G$ are irreducible and attain the minimal possible Gelfand--Kirillov dimension among all infinite-dimensional unitary representations. By Berger's classification \cite{Ber57} the symmetric subgroups of $G$ are
\begin{align*}
 K &= \upU(n),\\
 H_1 &= \GL(p,\CC)\times\GL(q,\CC), && n=p+q,\,p,q\geq1,\\
 H_2 &= \upU(p,q), && n=p+q,\,p,q\geq1,\\
 H_3 &= \Sp(m,\CC), && n=2m,\\
 H_4 &= \GL(m,\HH), && n=2m,\\
 H_5 &= \upO(n,\CC),\\
 H_6 &= \GL(n,\RR),
\end{align*}
where $H_3$ and $H_4$ only occur for even $n$. In this paper we give the explicit branching laws for the restriction of $\pi_{i\lambda,m}^G$ to these subgroups (see Theorems~\ref{thm:BranchingK} to \ref{thm:BranchingH6}).

The main machinery to derive the branching laws is Mackey theory. Following \cite{KOP11}, we denote for a homogeneous $G$-space $X$ and a homogeneous line bundle $\calL\to X$ by $L^2(X,\calL)$ the space of $L^2$-sections of the line bundle $\calL\otimes(\bigwedge^{\dim X}T^*X)^{\frac{1}{2}}$. Then the representations $\pi_{i\lambda,k}$ are realized on $L^2(G/P,\calL_{i\lambda,k})$ for $\calL_{i\lambda,k}=G\times_P\chi_{i\lambda,k}$. Each symmetric subgroup $H\subseteq G$ has by a general result of Wolf~\cite{Wol74} finitely many orbits on $G/P$ and hence there exist open orbits $\calO^1,\ldots,\calO^s$ such that their union is dense in $G/P$ (in all our cases $s=1$ or $s=2$). By Mackey theory restriction to these open orbits defines an $H$-equivariant unitary isomorphism
\[ L^2(G/P,\calL_{i\lambda,k}) \cong \bigoplus_{j=1}^s L^2(\calO^j,\calL_{i\lambda,k}|_{\calO^j}) \]
and hence the branching problem is equivalent to the decomposition of the unitary representations $L^2(\calO^j,\calL_{i\lambda,k}|_{\calO^j})$ into irreducible $H$-representations. In all our cases $\calO^j$ is a symmetric space or a flag variety or a fibration over one of those. This allows to apply existing results such as Plancherel formulas for reductive symmetric spaces or structure theory for parabolically induced representations.

We remark that the subgroups $H_1$, $H_2$, $H_3$ and $H_4$ are analogues of the four non-compact symmetric subgroups of $\GL(n,\RR)$ considered in \cite{KOP11}. The proofs of the corresponding branching laws are of the same nature as those in \cite{KOP11}. However, the two groups $H_5$ and $H_6$ do not have an analogue in $\GL(n,\RR)$. We still use Mackey theory for these cases, but here different phenomena occur. For $H_5$ we employ the Plancherel formula for the semisimple symmetric space $\upO(n,\CC)/(\upO(1,\CC)\times \upO(n-1,\CC))$ due to Delorme~\cite{Del98} and van den Ban--Schlichtkrull~\cite{vdBS05a,vdBS05b}. To obtain the branching law for $H_6$ we use the Plancherel formula for line bundles over the upper half plane $\SL(2,\RR)/\SO(2)$ due to Shimeno~\cite{Shi94}.

\subsection*{Acknowledgement.} We thank B. {\O}rsted for helpful discussions on the topic of this paper.
\section{Unitary principal series representations}\label{sec:UnitaryPrincipalSeriesRepresentations}

Let $G=\GL(n,\CC)$, $n\geq2$, and denote by $P$ the standard maximal parabolic subgroup of $G$ corresponding to the partition $(1,n-1)$ of $n$. Then $P\cong(\GL(1,\CC)\times\GL(n-1,\CC))\ltimes\CC^{n-1}$ and $G/P\cong\CC\PP^{n-1}$, the complex projective space of dimension $n-1$. We use the notation
\begin{equation*}
 [z] = [z_1:\dots:z_n] = \CC z\in\CC\PP^{n-1},
\end{equation*}
for $z\in\CC^n$, $z\neq0$. The complex projective space $\CC\PP^{n-1}$ is also homogeneous under the action of the maximal compact subgroup $K=\upU(n)$ and hence $\CC\PP^{n-1}\cong K/M$ with $M=K\cap P=\upU(1)\times \upU(n-1)$.

We define for $t\in\CC$ and $k\in\ZZ$ a character $\chi_{t,k}$ of $\GL(1,\CC)\times\GL(n-1,\CC)$ by
\begin{align*}
 \chi_{t,k}(\diag(w,g)) &:= |w|^t\left(\frac{w}{|w|}\right)^k, & w\in\GL(1,\CC),\ g\in\GL(n-1,\CC),
\end{align*}
and extend it to $P$ by letting the nilpotent radical act trivially. Let $\pi_{t,k}:=\Ind_P^G(\chi_{t,k})$ denote the corresponding unitary principal series representations (normalized parabolic induction). The structure of these representations was studied by Dooley--Zhang~\cite{DZ97} (see also \cite{HL99,MS12}). For $t=i\lambda\in i\RR$ and $k\in\ZZ$ the representation $\pi_{t,k}$ is unitary and irreducible and called \textit{unitary principal series representation}. If $n\geq3$ the unitary principal series are the only non-trivial unitary constituents of the principal series (see e.g. \cite[Corollary 2.4.3 \&\ 2.5.8]{HL99}). In \cite[Theorem 5.3]{MS12} we showed that the unitary principal series representations $\pi_{i\lambda,k}$ attain the minimal possible Gelfand--Kirillov dimension among all infinite-dimensional unitary irreducible representations of $G$. (In fact, we only showed the statement for $k=0$, but the proof carries over to the case $k\in\ZZ$.)

The unitary representation $\pi_{i\lambda,k}$ can be realized on $L^2(G/P,\calL_{i\lambda,k})$, where $\calL_{i\lambda,k}=G\times_P\chi_{i\lambda,k}$ is the line bundle associated to the character $\chi_{i\lambda,k}$ and $G$ acts via the left-regular representation. Considered as $K$-homogeneous bundle over $G/P\cong K/M$ the line bundle $\calL_{i\lambda,k}$ is the bundle $\calL_k$ corresponding to the character $(w,g)\mapsto w^k$ of $M=\upU(1)\times \upU(n-1)$. Hence as $K$-representations
\begin{equation*}
 L^2(G/P,\calL_{i\lambda,k})\cong L^2(K/M,\calL_k).
\end{equation*}
Since $K/M=\CC\PP^{n-1}$ we can identify $L^2(K/M,\calL_k)$ with the subspace of $L^2(\upU(n)/\upU(n-1))$ consisting of functions on the unit sphere $\upU(n)/\upU(n-1)\cong S^{2n-1}\subseteq\CC^n$ which are homogeneous of degree $-k$ under the action of $\upU(1)$ on $S^{2n-1}$. It is well-known that as $\upO(2n)$-representations the space $L^2(S^{2n-1})$ decomposes into the Hilbert space direct sum
\begin{align*}
 L^2(S^{2n-1}) &= {\sum_{j=0}^\infty}\raisebox{0.15cm}{$^\oplus$}{\,\calH^j(\RR^{2n})},
\end{align*}
where $\calH^j(\RR^{2n})$ denotes the irreducible representation of $\upO(2n)$ on the space of spherical harmonics on $\RR^{2n}$ of degree $j$. Under the action of $\upU(n)\subseteq\upO(2n)$ each space of spherical harmonics decomposes into $\upU(n)$-irreducibles as
\begin{align*}
 \calH^j(\RR^{2n}) &= \bigoplus_{\substack{\alpha,\beta\geq0\\\alpha+\beta=j}}{\calH^{\alpha,\beta}(\CC^n)},
\end{align*}
where $\calH^{\alpha,\beta}(\CC^n)$ denotes the subspace of $\calH^j(\RR^{2n})$ consisting of polynomials which are homogeneous of degree $\alpha$ in $z$ and homogeneous of degree $\beta$ in $\overline{z}$. This implies the following branching law:

\begin{theorem}\label{thm:BranchingK}
Upon restriction to the maximal compact subgroup $K=\upU(n)$ the representation $\pi_{i\lambda,k}^G$ of $G=\GL(n,\CC)$ decomposes as
\begin{align*}
 \pi_{i\lambda,k}^G|_K \cong {\sum_{\substack{\alpha,\beta\geq0\\\alpha-\beta=-k}}}\hspace{-.3cm}\raisebox{0.15cm}{$^\oplus$}{\,\calH^{\alpha,\beta}(\CC^n)}.
\end{align*}
\end{theorem}
\section{Branching law for $\GL(p+q,\CC)\searrow\GL(p,\CC)\times\GL(q,\CC)$}

The subgroup $H=H_1=\GL(p,\CC)\times\GL(q,\CC)$, $n=p+q$, has the open dense orbit
\begin{align*}
 \calO &= H\cdot[z_0] = \{[z':z'']:z'\in\CC^p\setminus\{0\},z''\in\CC^q\setminus\{0\}\}
\end{align*}
on $G/P=\CC\PP^{n-1}$ where
\[ z_0 = (\underbrace{1,0,\ldots,0}_{p},\underbrace{1,0,\ldots,0}_{q})\in\CC^{p+q}. \]
The stabilizer of $[z_0]$ in $H$ is
\begin{multline}
 S = \Bigg\{g=\left(\begin{array}{cccc}w&a&0&0\\0&g_1&0&0\\0&0&w&b\\0&0&0&g_2\end{array}\right):w\in\CC^\times,g_1\in\GL(p-1,\CC),\\
 g_2\in\GL(q-1,\CC),a\in\CC^{p-1},b\in\CC^{q-1}\Bigg\}.\label{eq:StabilizerGLGL}
\end{multline}
Restriction of the bundle $\calL_{i\lambda,k}$ to $\calO\cong H/S$ yields the line bundle induced by the character $S\to\CC^\times,\,g\mapsto|w|^{i\lambda}(\frac{w}{|w|})^k$ with $g\in S$ as in \eqref{eq:StabilizerGLGL}. The orbit $\calO$ has an $H$-equivariant fibration
\begin{align*}
 \calO\to\CC\PP^{p-1}\times\CC\PP^{q-1},\,[z':z'']\mapsto([z'],[z''])
\end{align*}
with fiber $\GL(1,\CC)=\CC^\times=\RR_+\times\SS^1$. Taking the Mellin transform along $\RR_+$ and the Fourier series expansion on $\SS^1$ we obtain the branching law. For the formulation we denote by $\pi_{i\lambda,k}^{\GL(p,\CC)}$ and $\pi_{i\lambda,k}^{\GL(q,\CC)}$ the corresponding unitary principal series representations of $\GL(p,\CC)$ and $\GL(q,\CC)$ as defined in Section~\ref{sec:UnitaryPrincipalSeriesRepresentations}.

\begin{theorem}
Upon restriction to the symmetric subgroup $H_1=\GL(p,\CC)\times\GL(q,\CC)$ the representation $\pi_{i\lambda,k}^G$ of $G=\GL(p+q,\CC)$ decomposes as
\begin{align*}
 \pi_{i\lambda,k}^G|_{H_1} &\cong {\sum_{k'\in\ZZ}^\infty}\raisebox{0.15cm}{$^\oplus$}\,{\int_\RR^\oplus{\pi_{i\lambda',k'}^{\GL(p,\CC)}\otimes\pi_{i(\lambda-\lambda'),k-k'}^{\GL(q,\CC)}\td\lambda'}}.
\end{align*}
\end{theorem}
\section{Branching law for $\GL(p+q,\CC)\searrow \upU(p,q)$}\label{sec:BranchingU(p,q)}

Let $H=H_2=\upU(p,q)$, $p,q\geq1$, $n=p+q$, be the subgroup of all elements in $\GL(n,\CC)$ which preserve the sesquilinear form
\begin{align*}
 (z|w)_{p,q} &:= z_1\overline{w_1}+\cdots+z_p\overline{w_p}-z_{p+1}\overline{w_{p+1}}-\cdots-z_{p+q}\overline{w_{p+q}}, & z,w\in\CC^{p+q}.
\end{align*}
The group $H$ has exactly two open orbits on $G/P=\CC\PP^{n-1}$, namely
\begin{align*}
 \calO^+ &:= H\cdot[1:0:\dots:0] = \{[z]:(z|z)_{p,q}>0\},\\
 \calO^- &:= H\cdot[0:\dots:0:1] = \{[z]:(z|z)_{p,q}<0\}.
\end{align*}
Both orbits are semisimple symmetric spaces, more precisely $\calO^\pm\cong H/S^\pm$ with
\begin{align*}
 S^+ &:= \upU(1)\times \upU(p-1,q), && \mbox{and} & S^- &:= \upU(p,q-1)\times \upU(1).
\end{align*}
The line bundle $\calL_{i\lambda,k}$ restricted to $\calO^+$ is the line bundle $\calL_k^+$ induced by the character $(z,g)\mapsto z^k$ of $S^+=\upU(1)\times \upU(p-1,q)$, and restricted to $\calO^-$ it is the line bundle $\calL_k^-$ induced by the character $(g,z)\mapsto z^k$ of $S^-=\upU(p,q-1)\times \upU(1)$. Hence, as $H$-representations we have
\begin{align*}
 L^2(G/P,\calL_{i\lambda,k}) &\cong L^2(H/S^+,\calL_k^+)\oplus L^2(H/S^-,\calL_k^-).
\end{align*}
To describe the decomposition of $L^2(H/S^\pm,\calL_k^\pm)$ into irreducible unitary $H$-representations we use the degenerate principal series of $H$. Let $P_H=M_HA_HN_H$ be the parabolic subgroup of $H$ with
\begin{align*}
 M_H &= \left\{\diag(\mu,g,\mu):\mu\in \upU(1),g\in \upU(p-1,q-1)\right\},\\
 A_H &= \exp\RR A_0,\\
 N_H &= \exp(\frakh_{\alpha}+\frakh_{2\alpha}),
\end{align*}
where
\begin{align*}
 A_0 &= \left(\begin{array}{ccc}&&1\\&&\\1&&\end{array}\right)
\end{align*}
and $\alpha\in\fraka_\CC^*$ is defined by $\alpha(A_0):=1$. For $t\in\CC$ and $k\in\ZZ$ let $\omega_{t,k}$ denote the character of $P_H=M_HA_HN_H$ given by
\begin{align*}
 \omega_{t,k}(\diag(\mu,g,\mu)e^{sA_0}n) &:= \mu^k e^{st}.
\end{align*}
We consider the induced representation $\pi_{t,k}^{\upU(p,q)}:=\Ind_{P_H}^H(\omega_{t,k})$. For $t\in i\RR$ this representation belongs to the unitary principal series and is irreducible except possibly for $t=0$ and $k\in n+1+2\ZZ$ (where it might decompose into two components, see \cite[Section 6]{MvD03} and \cite[Section 4.5]{HT93}). To describe the unitary subrepresentations of $\pi_{t,k}^{\upU(p,q)}$ that occur discretely in the decomposition of $L^2(H/S^\pm,\calL_k^\pm)$ we distinguish three cases:

\begin{enumerate}
\item Assume $p,q>1$. Then for $t<0$ and $t\in k+n+1+2\ZZ$ the representation $\pi_{t,k}^{\upU(p,q)}$ has two irreducible unitary subquotients (see e.g. \cite[Section 6]{MvD03}). Let $\pi_{t,k,\pm}^{\upU(p,q)}$ denote the irreducible subquotient with $\upU(p)\times \upU(q)$-type decomposition
\begin{equation}
 \pi_{t,k,\pm}|_{\upU(p)\times \upU(q)} =  {\sum_{\substack{j,\ell\in\NN_0\\\pm(j-\ell+p-q)>|t|}}}\hspace{-8mm}\raisebox{0.15cm}{$^\oplus$}\hspace{7mm}{\, \bigoplus_{\substack{0\leq\alpha\leq j,0\leq\beta\leq\ell\\2(\alpha+\beta)=j+\ell-k}}{\,\calH^{\alpha,j-\alpha}(\CC^p)\otimes\calH^{\beta,\ell-\beta}(\CC^q)}}.\label{eq:KtypesUpqBigConstituents}
\end{equation}
We remark that identifying $\upU(q,p)\cong \upU(p,q)$ we have $\pi_{t,k,-}^{\upU(p,q)}\cong\pi_{t,k,+}^{\upU(q,p)}$.

\item Assume $p>1$ and $q=1$. Since $\calH^{\beta_1,\beta_2}(\CC)=0$ for $\beta_1\beta_2\neq0$ we put
\[ \calH^\beta(\CC) = \begin{cases}\calH^{\beta,0}(\CC) & \mbox{for $\beta\geq0$,}\\\calH^{0,-\beta}(\CC) & \mbox{for $\beta\leq0$.}\end{cases} \]
For $|k|<p$ there is a complementary series, more precisely $\pi_{t,k}^{\upU(p,1)}$ is unitary and irreducible for $|t|<p-|k|$ (see e.g. \cite[Section 4.5]{HT93}). To simplify notation we denote those representations also by $\pi_{t,k,+}^{\upU(p,1)}=\pi_{t,k}^{\upU(p,1)}$.

For $t<0$ and $t\in k+n+1+2\ZZ$ not in the range of the possible complementary series (i.e. $t<-(p-|k|)$) there is always a ``big'' unitarizable subquotient in $\pi_{t,k}^{\upU(p,1)}$ denoted by $\pi_{t,k,+}^{\upU(p,1)}$ with $K$-types given by the formula \eqref{eq:KtypesUpqBigConstituents}. For $t>-p-|k|$ this is the unique irreducible constituent containing the $K$-type $\calH^{0,k}(\CC^p)\otimes\calH^0(\CC)$ for $k\geq0$ and $\calH^{-k,0}(\CC^p)\otimes\calH^0(\CC)$ for $k\leq0$.

For $t\in k+n+1+2\ZZ$ with $p-|k|\leq t<0$ there is further a small unitarizable subquotient of $\pi_{t,k}^{\upU(p,1)}$ denoted by $\pi_{t,k,-}^{\upU(p,1)}$ with $K$-types given by
\[ \pi_{t,k,-}^{\upU(p,1)}|_{\upU(p)\times \upU(1)} = \begin{cases}{\displaystyle\sum_{\substack{\alpha_1,\alpha_2\in\NN_0,\beta\in\ZZ\\\alpha_1+\alpha_2+\beta\leq t-p\\\alpha_1-\alpha_2+\beta=-k}}}
\hspace{-7mm}\raisebox{0.15cm}{$^\oplus$}\,\calH^{\alpha_1,\alpha_2}(\CC^p)\otimes\calH^\beta(\CC) & \mbox{for $k>0$,}\\[11mm]\displaystyle
{\sum_{\substack{\alpha_1,\alpha_2\in\NN_0,\beta\in\ZZ\\\alpha_1+\alpha_2-\beta\leq t-p\\\alpha_1-\alpha_2+\beta=-k}}}\hspace{-7mm}\raisebox{0.15cm}{$^\oplus$}\,\calH^{\alpha_1,\alpha_2}(\CC^p)\otimes\calH^\beta(\CC) & \mbox{for $k<0$.}\end{cases} \]
Note that this only occurs when $|k|>p$. In this case $\pi_{t,k,-}^{\upU(p,1)}$ is the unique irreducible constituent containing the $K$-type $\calH^{0,0}(\CC^p)\otimes\calH^{-k}(\CC)$. For $k>0$ the representations $\pi_{t,k,-}^{\upU(p,1)}$ belong to the holomorphic discrete series and for $k<0$ to the antiholomorphic discrete series.

\item Assume $p=1$ and $q>1$. Identifying $\upU(1,q)\cong \upU(q,1)$ we use (2) to define $\pi_{t,k,+}^{\upU(1,q)}:=\pi_{t,k,-}^{\upU(q,1)}$ and $\pi_{t,k,-}^{\upU(1,q)}:=\pi_{t,k,+}^{\upU(q,1)}$.
\end{enumerate}

Let
\[ A_+^k(p,q) := \begin{cases}(k+n+1+2\ZZ)\cap(-\infty,0) & \mbox{for $p>1$, $q\geq1$,}\\(k+n+1+2\ZZ)\cap(-(|k|-q),0) & \mbox{for $p=1$, $q\geq1$,}\end{cases} \]
and put $A_-^k(p,q):=A_+^k(q,p)$. Here we use the convention $(x,0)=\emptyset$ for $x\geq0$. Hence $A_+^k(1,q)=A_-^k(q,1)=\emptyset$ for $|k|\leq q$. Then by \cite[Theorem 10.3]{MvD03} the space $L^2(H/S^\pm,\calL_k^\pm)$ decomposes into
\begin{align*}
 L^2(H/S^\pm,\calL_k^\pm) \cong {\sum_{t\in A_\pm^k(p,q)}}
\hspace{-3.8mm}\raisebox{0.15cm}{$^\oplus$}\,{\pi_{t,k,\pm}^{\upU(p,q)}}\oplus\int_{i\RR_+}^\oplus{\pi_{t,k}^{\upU(p,q)}\td t}.
\end{align*}
This gives the following result:

\begin{theorem}
Upon restriction to the symmetric subgroup $H_2=\upU(p,q)$ the representation $\pi_{i\lambda,k}^G$ of $G=\GL(p+q,\CC)$ decomposes as
\begin{align*}
 \pi_{i\lambda,k}^G|_{H_2} &\cong {\sum_{t\in A_+^k(p,q)}}
\hspace{-3.8mm}\raisebox{0.15cm}{$^\oplus$}\,{\pi_{t,k,+}^{\upU(p,q)}} \oplus {\sum_{t\in A_-^k(p,q)}}
\hspace{-3.8mm}\raisebox{0.15cm}{$^\oplus$}\,{\pi_{t,k,-}^{\upU(p,q)}} \oplus 2\!\int_{i\RR_+}^\oplus{\pi_{t,k}^{\upU(p,q)}\td t}.
\end{align*}
\end{theorem}

We remark that for $n>2$ there always occur infinitely many discrete components in the branching law. For $\min(p,q)=1$ there are complementary series representations among them if $|k|<\max(p,q)$ and holomorphic or antiholomorphic discrete series representations if $|k|>\max(p,q)$ (finitely many in both cases).
\section{Branching law for $\GL(2m,\CC)\searrow\Sp(m,\CC)$}

For $n=2m$ the subgroup $H=H_3=\Sp(m,\CC)$ of $G$ acts transitively on $G/P\cong\CC\PP^{n-1}$ and $G/P\cong H/P_H$, where $P_H=H\cap P=L_HN_H$ is a maximal parabolic subgroup of $H$ with
\begin{align}
 L_H &= \left\{\ell=\left(\begin{array}{cccc}w&0&0&0\\0&A&0&B\\0&0&w^{-1}&0\\0&C&0&D\end{array}\right):w\in\GL(1,\CC),\left(\begin{array}{cc}A&B\\C&D\end{array}\right)\in\Sp(m-1,\CC)\right\}\label{eq:SymplecticLH}
\end{align}
and $N_H$ a complex Heisenberg group of dimension $2m-1$. The line bundle $\calL_{i\lambda,k}$ over $G/P\cong H/P_H$ is, as an $H$-homogeneous bundle, induced from the character $\omega_{i\lambda,k}$ of $P_H$ given by
\begin{align*}
 \omega_{i\lambda,k}(\ell n) &= |w|^{i\lambda}(\tfrac{w}{|w|})^k
\end{align*}
for $\ell\in L_H$ as in \eqref{eq:SymplecticLH} and $n\in N_H$. Hence $\pi_{i\lambda,k}^G|_H\cong\Ind_{P_H}^H(\omega_{i\lambda,k})$. By \cite[Theorem 7]{Gro71} the unitary representations $\pi_{i\lambda,k}^H:=\Ind_{P_H}^H(\omega_{i\lambda,k})$ are irreducible for $(i\lambda,k)\neq(0,0)$ and decompose into two irreducible components for $(i\lambda,k)=(0,0)$. In \cite[Theorem 1]{Cla12} these two components are characterized in terms of their $K$-types:
\begin{align*}
 \pi_{0,0}^{\Sp(m,\CC)} &= \pi_{0,0,+}^{\Sp(m,\CC)}\oplus\pi_{0,0,-}^{\Sp(m,\CC)},
\end{align*}
where
\[ \pi_{0,0,+}^{\Sp(m,\CC)}|_{\Sp(m)} \cong {\sum_{\substack{\alpha-\beta\geq0\\\alpha-\beta\in4\ZZ}}}
\hspace{-2.5mm}\raisebox{0.15cm}{$^\oplus$}\,{\calH^{\alpha,\beta}(\HH^m)}, \qquad 
\pi_{0,0,-}^{\Sp(m,\CC)}|_{\Sp(m)} \cong {\sum_{\substack{\alpha-\beta\geq0\\\alpha-\beta\in2+4\ZZ}}}
\hspace{-4.45mm}\raisebox{0.15cm}{$^\oplus$}\,{\calH^{\alpha,\beta}(\HH^m)}. 
\]
Here $\calH^{\alpha,\beta}(\HH^m)$ denotes the irreducible representation of $\Sp(m)$ of highest weight $(\alpha,\beta,0,\ldots,0)$ in the standard notation. It occurs with multiplicity $\alpha-\beta+1$ in the decomposition of $\calH^{\alpha+\beta}(\RR^{4m})$ into irreducible $\Sp(m)$-representations.

\begin{theorem}
Upon restriction to the symmetric subgroup $H_3=\Sp(m,\CC)$ the representation $\pi_{i\lambda,k}^G$ of $G=\GL(2m,\CC)$ is
\begin{itemize}
\item for $(i\lambda,k)\neq(0,0)$ irreducible and isomorphic to $\pi_{i\lambda,k}^{\Sp(m,\CC)}$,
\item for $(i\lambda,k)=(0,0)$ reducible and decomposes into two irreducible components:
\begin{align*}
 \pi_{0,0}^G|_{H_3} &\cong \pi_{0,0,+}^{\Sp(m,\CC)}\oplus\pi_{0,0,-}^{\Sp(m,\CC)}.
\end{align*}
\end{itemize}
\end{theorem}
\section{Branching law for $\GL(2m,\CC)\searrow\GL(m,\HH)$}

For $n=2m$ the group $H=H_4=\GL(m,\HH)\subseteq\GL(2m,\CC)=G$ acts transitively on $\CC^{2m}\setminus\{0\}$ and hence also transitively on $G/P\cong\CC\PP^{n-1}$. Thus $G/P$ is identified with the homogeneous space $H/(H\cap P)$, where
\begin{align*}
 H\cap P &= \left\{\left(\begin{array}{cc}a&b\\0&C\end{array}\right):a\in\CC^\times,b\in\HH^{m-1},C\in\GL(m-1,\HH)\right\},
\end{align*}
viewed as $m\times m$ matrices over the quaternions $\HH$. Let $P_H$ be the maximal parabolic subgroup of $H$ defined by
\begin{align*}
 P_H &:= \left\{\left(\begin{array}{cc}a&b\\0&C\end{array}\right):a\in\HH^\times,b\in\HH^{m-1},C\in\GL(m-1,\HH)\right\}.
\end{align*}
Then $H/(H\cap P)\to H/P_H$ is an $H$-equivariant fibration with fiber $\HH^\times/\CC^\times\cong\Sp(1)/\upU(1)$, where we identify $\Sp(1)$ with the group of unit quaternions in $\HH$ and $\upU(1)$ with the group of complex numbers of absolute value one. As an $H$-homogeneous bundle the line bundle $\calL_{i\lambda,k}$ is induced from the character $\chi_{i\lambda,k}|_{H\cap P}$. Using induction in stages we find
\begin{equation*}
 \pi_{i\lambda,k}^G|_H = \Ind_{H\cap P}^H(\chi_{i\lambda,k}|_{H\cap P}) = \Ind_{P_H}^H(\Ind_{H\cap P}^{P_H}(\chi_{i\lambda,k}|_{H\cap P})).
\end{equation*}
Let $\calL_k$ be the line bundle over $\Sp(1)/\upU(1)$ induced by the character $z\mapsto z^k$ of $\upU(1)$. Then
\begin{align*}
 L^2(\Sp(1)/\upU(1),\calL_k) \cong {\sum_{\substack{j\geq|k|\\j\in k+2\ZZ}}}
\hspace{-2.2mm}\raisebox{0.15cm}{$^\oplus$}\,{V_j},
\end{align*}
where $(\tau_j,V_j)$ is the (unique) irreducible representation of $\Sp(1)\cong\SU(2)$ of dimension $j+1$. We define a unitary representation $\omega_{i\lambda,j}$ of $P_H$ on $V_j$ by
\begin{equation*}
 \omega_{i\lambda,j}\left(\begin{array}{cc}a&b\\0&C\end{array}\right) := |a|^{i\lambda}\tau_j(\tfrac{a}{|a|})
\end{equation*}
and let $\pi_{i\lambda,j}^H:=\Ind_{P_H}^H(\omega_{i\lambda,j})$. By \cite[Theorems 3.3 \&\ 4.3]{Pas99} the representations $\pi_{i\lambda,j}^H$ are irreducible unitary representations of $H$. Hence we obtain the decomposition of $\pi_{i\lambda,k}^G$ into irreducible $H$-representations:

\begin{theorem}
Upon restriction to the symmetric subgroup $H_4=\GL(m,\HH)$ the representation $\pi_{i\lambda,k}^G$ of $G=\GL(2m,\CC)$ decomposes as
\begin{align*}
 \pi_{i\lambda,k}^G|_{H_4} &\cong {\sum_{\substack{j\geq|k|\\j\in k+2\ZZ}}}
\hspace{-2.2mm}\raisebox{0.15cm}{$^\oplus$}\,{\pi_{i\lambda,j}^{\GL(m,\HH)}}.
\end{align*}
\end{theorem}
\section{Branching law for $\GL(n,\CC)\searrow \upO(n,\CC)$}

The group $H=H_5=\upO(n,\CC)$ has the open dense orbit
\begin{align*}
 \calO &= H\cdot[1:0:\dots:0] = \{[z]:z\in\CC^n,\ z_1^2+\cdots+z_n^2\neq0\}
\end{align*}
which is the homogeneous space $\calO\cong H/S$ with $S=H\cap P=\upO(1,\CC)\times \upO(n-1,\CC)$. The restriction of the line bundle $\calL_{i\lambda,k}$ to $\calO$ is the $H$-equivariant line bundle $\calL_{k+2\ZZ}$ where for $\delta\in\ZZ/2\ZZ$ we denote by $\calL_\delta$ the $H$-equivariant line bundle over $H/S$ associated to the character $(z,g)\mapsto z^\delta$ of $S$. Hence
\begin{align*}
 L^2(G/P,\calL_{i\lambda,k}) \cong L^2(H/S,\calL_{k+2\ZZ}).
\end{align*}
Put
\begin{equation*}
 S_1 := S = \upO(1,\CC)\times \upO(n-1,\CC), \qquad S_2 := \{1\}\times \upO(n-1,\CC),
\end{equation*}
then
\begin{equation*}
 L^2(H/S_1) \cong L^2(H/S,\calL_0), \qquad L^2(H/S_2) \cong L^2(H/S,\calL_0)\oplus L^2(H/S,\calL_1).
\end{equation*}
Therefore it suffices to find the Plancherel formula for $L^2(H/S_1)$ and $L^2(H/S_2)$. Both $H/S_1$ and $H/S_2$ are semisimple symmetric spaces and we can use the Plancherel formula by Delorme~\cite{Del98} and van den Ban--Schlichtkrull~\cite{vdBS05a,vdBS05b}. We follow the outline in \cite{vdB05}.

We fix the Cartan involution $\theta(g):=(g^*)^{-1}=\overline{g}$ of $H$ with corresponding maximal compact subgroup $H^\theta=H\cap K=\upO(n)$. The involution $\sigma$ of $H$ given by $\sigma(g):=\1_{1,n-1}\circ g\circ\1_{1,n-1}$ with $\1_{1,n-1}:=\diag(1,-\1_{n-1})$ satisfies $H^\sigma_0\subseteq S_1,S_2\subseteq H^\sigma$ and hence $S_1$ and $S_2$ are symmetric subgroups of $H$. We note that for $i=1,2$ we have $\rk(H/S_i)=2$ whereas $\rk((H\cap K)/(S_i\cap K))=1$ and hence there occurs no discrete spectrum in $L^2(H/S_i)$ by a result of Flensted-Jensen \cite{FJ80} and Oshima--Matsuki \cite{OM84}. To find the continuous spectrum we observe that there is (up to Weyl group action) only one non-trivial $\sigma\theta$-stable parabolic subgroup $P_H=M_HA_HN_H$ of $H$. It has abelian nilradical $N_H\cong\CC^{n-2}$ and
\begin{align*}
 M_H &= \upO(2)\times \upO(n-2,\CC), & A_H &= \exp(\fraka_H)
\end{align*}
with
\begin{align*}
 \fraka_H &:= \RR A_0, & A_0 &:= \diag(\left(\begin{array}{cc}0&i\\-i&0\end{array}\right),0,\ldots,0).
\end{align*}
Note that $\fraka_H$ is already maximal abelian in $\{X\in\frakh:\theta(X)=\sigma(X)=-X\}$. Let $W=N_{H\cap K}(\fraka_H)/Z_{H \cap K}(\fraka_H)=(\upO(2)\times \upO(n-2))/(\SO(2)\times \upO(n-2))=\{\1,w_0\}$ with $w_0=\diag(-1,1,\ldots,1)$. For $i=1,2$ we denote by $W_i$ the natural image of the group $N_{S_i\cap K}(\fraka_H)$ in $W$. It is easy to see that $W_i=W$ for both $i=1,2$. Therefore, following \cite[Section 8]{vdB05}, we have to consider discrete series for the spaces $X_i=M_H/(S_i\cap M_H)$ and their $(S_i\cap M_H)$-fixed vectors. We find
\begin{align*}
 X_1 &\cong \upO(2)/(\upO(1)\times\upO(1)), & X_2 &\cong \upO(2)/(\{1\}\times\upO(1)).
\end{align*}
For $i=1$ the $(\upO(1)\times\upO(1))$-spherical representations of $\upO(2)$ are the spherical harmonics $\calH^j(\RR^2)$ of even degree $j\in2\NN_0$. For $i=2$ the $(\{1\}\times\upO(1))$-spherical representations of $\upO(2)$ are spherical harmonics $\calH^j(\RR^2)$ of arbitrary degree $j\in\NN_0$. Hence
\[ L^2(X_1) = {\sum_{j\in2\NN_0}}
\hspace{-1.2mm}\raisebox{0.15cm}{$^\oplus$}\,\calH^j(\RR^2), \qquad L^2(X_2) = {\sum_{j\in\NN_0}}
\hspace{-.5mm}\raisebox{0.15cm}{$^\oplus$}\,\calH^j(\RR^2). \]
Note that $\dim\calH^0(\RR^2)=1$ and $\dim\calH^j(\RR^2)=2$ for $j>0$. Extend the representations $\calH^j(\RR^2)$ of $\upO(2)$ to $M_H=\upO(2)\times \upO(n-2,\CC)$ by letting $\upO(n-2,\CC)$ act trivially. Further, identify $(\fraka_H)_\CC^*$ with $\CC$ by $\alpha\mapsto\alpha(A_0)$. The induced representations
\begin{equation*}
 \pi_{t,j}^H := \Ind_{P_H}^H(\calH^j(\RR^2)\otimes e^t\otimes\1)
\end{equation*}
are unitary irreducible and pairwise inequivalent for $t\in i\RR_+$ and $j\in\NN_0$. The Plancherel formula \cite[Theorem 10.15]{vdB05} gives
\begin{align*}
 L^2(H/S_1) &\cong {\sum_{j\in2\NN_0}}
\hspace{-1.2mm}\raisebox{0.15cm}{$^\oplus$}\,{\int^\oplus_{i\RR_+}{\pi_{t,j}^H\td t}}, & L^2(H/S_2) &\cong {\sum_{j\in\NN_0}}
\hspace{-.5mm}\raisebox{0.15cm}{$^\oplus$}\,{\int^\oplus_{i\RR_+}{\pi_{t,j}^H\td t}}.
\end{align*}
Putting everything together gives:

\begin{theorem}
Upon restriction to the symmetric subgroup $H_5=\upO(n,\CC)$ the representation $\pi_{i\lambda,k}^G$ of $G=\GL(n,\CC)$ decomposes as
\begin{align*}
 \pi_{i\lambda,k}^G|_{H_5} &\cong {\sum_{\substack{j\in\NN_0\\j\in k+2\ZZ}}}
\hspace{-2.2mm}\raisebox{0.15cm}{$^\oplus$}\,{\int_{i\RR_+}^\oplus{\pi_{t,j}^{\upO(n,\CC)}\td t}}.
\end{align*}
\end{theorem}
\section{Branching law for $\GL(n,\CC)\searrow\GL(n,\RR)$}

Let $H=H_6=\GL(n,\RR)$ and assume $n\geq3$. (The branching for $H=\GL(2,\RR)\cong U(1,1)$ was already treated in Section~\ref{sec:BranchingU(p,q)}.) The group $H$ has the open dense orbit
\begin{align*}
 \calO := H\cdot[1:i:0:\ldots:0] &= \{[ge_1+ige_2]:g\in\GL(n,\RR)\}\\
 &= \{[z]\in\CC\PP^{n-1}:\Re(z),\Im(z)\mbox{ linear independent}\}.
\end{align*}
on $G/P=\CC\PP^{n-1}$. The stabilizer subgroup of $[1:i:0:\dots:0]$ in $H$ is given by $S=(\RR_+\SO(2)\times\GL(n-2,\RR))\ltimes\RR^{2\times(n-2)}$. Let $\calL_{i\lambda,k}^H$ be the restriction of the line bundle $\calL_{i\lambda,k}$ to $\calO=H/S$ and $\omega_{i\lambda,k}$ be the associated character of $S$. Then $\omega_{i\lambda,k}$ is the trivial extension to $S$ of the character of $\RR_+\SO(2)$ given by
\begin{equation*}
 \omega_{i\lambda,k}\left(\begin{array}{cc}a&b\\-b&a\end{array}\right)=|a+ib|^{i\lambda}\left(\frac{a+ib}{|a+ib|}\right)^k.
\end{equation*}
Let $P_H^1:=(\GL(2,\RR)\times\GL(n-2,\RR))\ltimes\RR^{2\times(n-2)}$ then $H/S\to H/P_H^1$ is a fibration with fiber
\begin{equation*}
 \GL(2,\RR)/\RR_+\SO(2)\cong\{z\in\CC:\Im(z)\neq0\}.
\end{equation*}
We use induction in stages to find
\[ \pi_{i\lambda,k}^G|_H = L^2(H/S,\calL_{i\lambda,k}^H) = \Ind_S^H(\omega_{i\lambda,k}) = \Ind_{P_H^1}^H(\Ind_S^{P_H^1}(\omega_{i\lambda,k})). \]
Let $P_H^1=M_H^1A_H^1N_H^1$ be the Langlands decomposition of $P_H^1$. Then
\begin{align*}
 M_H^1 &= \left(\begin{array}{cc}\SLpm(2,\RR)&0\\0&\SLpm(n-2,\RR)\end{array}\right),\\
 A_H^1 &= \left(\begin{array}{cc}\RR_+\1_2&0\\0&\RR_+\1_{n-2}\end{array}\right),\\
 N_H^1 &= \left(\begin{array}{cc}\1_2&*\\0&\1_{n-2}\end{array}\right),
\end{align*}
where
\[ \SLpm(N,\RR) = \{g\in\GL(N,\RR):|\det(g)|=1\}. \]
Further we identify $(\fraka_H^1)_\CC^*$ with $\CC^2$ by
\[ (\fraka_H^1)_\CC^*\to\CC^2,\,\alpha\mapsto(\alpha(\diag(1,1,0,\ldots,0)),\alpha(\diag(0,0,1,\ldots,1))). \]
Then
\[ \pi_{i\lambda,k}^G|_H = \Ind_{P_H^1}^H(\Ind_{\SO(2)}^{\SLpm(2,\RR)}(e^{ik\theta})\otimes e^{(i\lambda,0)}\otimes\1), \]
where we view representations of $\SLpm(2,\RR)$ as representations of $M_H^1$ by letting the second factor $\SLpm(n-2,\RR)$ act trivially. We now find the decomposition of $\Ind_{\SO(2)}^{\SLpm(2,\RR)}(e^{ik\theta})$ into irreducible $\SLpm(2,\RR)$-representations.

Let $Q=M_QA_QN_Q\subseteq\SL(2,\RR)$ with
\[ M_Q = \left\{\pm\left(\begin{array}{cc}1&0\\0&1\end{array}\right)\right\}, \quad A_Q = \left\{\left(\begin{array}{cc}a&0\\0&a^{-1}\end{array}\right):a>0\right\}, \quad N_Q = \left\{\left(\begin{array}{cc}1&x\\0&1\end{array}\right):x\in\RR\right\}. \]
For $\varepsilon\in\ZZ/2\ZZ$ and $t\in\CC$ define a character $\varpi_{t,\varepsilon}$ of $Q$ by
\[ \varpi_{t,\varepsilon}\left(\begin{array}{cc}a&x\\0&a^{-1}\end{array}\right) := \sgn(a)^\varepsilon|a|^t, \qquad a\in\RR^\times,x\in\RR. \]
Let $\pi_{t,\varepsilon}^{\SL(2,\RR)}:=\Ind_Q^{\SL(2,\RR)}(\varpi_{t,\varepsilon})$ denote the corresponding induced representations. For $t\in i\RR_+$ the representation $\pi_{t,\varepsilon}^{\SL(2,\RR)}$ is unitary irreducible and for $t>0$ with $t\in\varepsilon+1+2\ZZ$ it contains two discrete series representations $\pi_{t,\varepsilon,\pm}^{\SL(2,\RR)}$ with $\SO(2)$-types $e^{\pm ij\theta}$, $j>t$. Then the Plancherel formula for line bundles over the upper half plane by Shimeno~\cite{Shi94} states that
\begin{align*}
 \Ind_{\SO(2)}^{\SL(2,\RR)}(e^{ik\theta}) = \bigoplus_{\substack{t>0\\t\in|k|-1-2\NN_0}}\pi_{t,k+2\ZZ,\sgn(k)}^{\SL(2,\RR)}\oplus\int_{i\RR_+}\pi_{t,k+2\ZZ}^{\SL(2,\RR)}\td t.
\end{align*}
Now let $Q^*=M_Q^*A_QN_Q$ be the parabolic subgroup of $\SLpm(2,\RR)$ with
\[ M_Q^* = \left\{\left(\begin{array}{cc}\pm1&0\\0&\pm1\end{array}\right)\right\}. \]
For $\varepsilon,\delta\in\ZZ/2\ZZ$ and $t\in\CC$ denote by $\varpi_{t,\varepsilon,\delta}$ the character of $Q^*$ given by
\[ \varpi_{t,\varepsilon,\delta}\left(\begin{array}{cc}a&x\\0&b\end{array}\right)\mapsto\sgn(a)^\varepsilon\sgn(b)^\delta|a|^t, \qquad |ab|=1, x\in\RR. \]
Denote by $\pi_{t,\varepsilon,\delta}^{\SLpm(2,\RR)}=\Ind_{Q^*}^{\SLpm(2,\RR)}(\varpi_{t,\varepsilon,\delta})$ the corresponding induced representation. Then it is easy to check that
\[ \pi_{t,\varepsilon,\delta}^{\SLpm(2,\RR)}\cong\pi_{t,\varepsilon+1,\delta+1}^{\SLpm(2,\RR)}, \]
whence we will write $\pi_{t,\varepsilon}^{\SLpm(2,\RR)}:=\pi_{t,\varepsilon,0}^{\SLpm(2,\RR)}$. This implies
\[ \Ind_{\SL(2,\RR)}^{\SLpm(2,\RR)}(\pi_{t,\varepsilon}^{\SL(2,\RR)}) \cong \bigoplus_{\delta=0,1}\pi_{t,\varepsilon+\delta,\delta}^{\SLpm(2,\RR)} \cong 2\cdot\pi_{t,\varepsilon}^{\SLpm(2,\RR)}. \]
Further, for $t>0$ with $t\in\varepsilon+1+2\ZZ$ the representation $\pi_{t,\varepsilon}^{\SLpm(2,\RR)}$ has a unique irreducible unitarizable subrepresentation $\pi_{t,\varepsilon,\ds}^{\SLpm(2,\RR)}$. It has $\upO(2)$-types $\calH^j(\RR^2)$, $j>t$, and we have
\[ \Ind_{\SL(2,\RR)}^{\SLpm(2,\RR)}(\pi_{t,\varepsilon,\pm}^{\SL(2,\RR)}) = \pi_{t,\varepsilon,\ds}^{\SLpm(2,\RR)}. \]
Putting things together we find that
\begin{align*}
 \Ind_{\SO(2)}^{\SLpm(2,\RR)}(e^{ik\theta}) &= \Ind_{\SL(2,\RR)}^{\SLpm(2,\RR)}(\Ind_{\SO(2)}^{\SL(2,\RR)}(e^{ik\theta}))\\
 &= \bigoplus_{\substack{t>0\\t\in|k|-1-2\NN_0}}\pi_{t,k+2\ZZ,\ds}^{\SLpm(2,\RR)}\oplus2\int_{i\RR_+}\pi_{t,k+2\ZZ}^{\SLpm(2,\RR)}\td t,
\end{align*}

Finally we note that by Harish-Chandra's Irreducibility Theorem (see e.g. \cite[Theorem 4.11]{KV96}) for $\varepsilon\in\ZZ/2\ZZ$ and $t>0$ with $t\in\varepsilon+1+2\ZZ$ we have
\[ \Ind_{P_H^1}^H(\pi_{t,\varepsilon,\ds}^{\SLpm(2,\RR)}\otimes e^{(i\lambda,0)}\otimes\1) \mbox{ is irreducible for $\lambda\neq0$.} \]
For the continuous spectrum we use \cite[Chapter VII, \S 2, Section 4]{Kna86} to find
\[ \Ind_{P_H^1}^H(\pi_{t,\varepsilon}^{\SLpm(2,\RR)}\otimes e^{(i\lambda,0)}\otimes\1) = \Ind_{P_H^2}^H(\xi_\varepsilon\otimes e^{(\frac{i\lambda+t}{2},\frac{i\lambda-t}{2},0)}\otimes\1), \]
where $P_H^2=M_H^2A_H^2N_H^2$ with
\begin{align*}
 M_H^2 &= \left(\begin{array}{ccc}\pm1&0&0\\0&\pm1&0\\0&0&\SLpm(n-2,\RR)\end{array}\right),\\
 A_H^2 &= \left(\begin{array}{ccc}\RR_+&0&0\\0&\RR_+&0\\0&0&\RR_+\1_{n-2}\end{array}\right),\\
 N_H^2 &= \left(\begin{array}{ccc}1&*&*\\0&1&*\\0&0&\1_{n-2}\end{array}\right),
\end{align*}
$\xi_\varepsilon$ is given by
\[ \xi_\varepsilon(\diag(\delta_1,\delta_2,g)) = \delta_1^\varepsilon, \]
and $(\fraka_H^2)_\CC^*$ is identified with $\CC^3$ by
\[ (\fraka_H^2)_\CC^*\to\CC^3, \quad \alpha\mapsto(\alpha(\diag(1,0,\ldots,0)),\alpha(\diag(0,1,0,\ldots,0)),\alpha(\diag(0,0,1,\ldots,1))). \]
Then again by Harish-Chandra's Irreducibility Theorem we find
\[ \Ind_{P_H^2}^H(\xi_\varepsilon\otimes e^{(\frac{i\lambda+t}{2},\frac{i\lambda-t}{2},0)}\otimes\1) \mbox{ is irreducible for $t\neq0,\pm i\lambda$,} \]
i.e. it is irreducible for generic $t\in i\RR$.

\begin{theorem}\label{thm:BranchingH6}
Upon restriction to the symmetric subgroup $H_6=\GL(n,\RR)$ the representation $\pi_{i\lambda,k}^G$ of $G=\GL(n,\CC)$ decomposes as
\begin{multline*}
 \pi_{i\lambda,k}^G|_{H_6} = \bigoplus_{\substack{t>0\\t\in|k|-1-2\NN_0}}\Ind_{P_H^1}^{\GL(n,\RR)}(\pi_{t,k+2\ZZ,\ds}^{M_H^1}\otimes e^{(i\lambda,0)}\otimes\1)\\
 \oplus \,2\int_{i\RR_+}\Ind_{P_H^2}^{\GL(n,\RR)}(\chi_{k+2\ZZ}\otimes e^{(i\lambda+t,i\lambda-t,0)}\otimes\1)\td t.
\end{multline*}
\end{theorem}

\bibliographystyle{amsplain}
\bibliography{bibdb}

\end{document}